\documentclass{amsart}
\usepackage{times,array,amssymb,tikz,mathrsfs,xcolor}
\usetikzlibrary{petri}

\newcounter{lemma}[section]

\newtheorem{Theorem}[lemma]{Theorem}

\newtheorem{Corollary}[lemma]{Corollary}
\newtheorem{Proposition}[lemma]{Proposition}

\newtheorem{theorem}{Theorem}

\theoremstyle{definition}

\newtheorem{Definition}[lemma]{Definition}

\newtheorem{Remark}[lemma]{Remark}

\numberwithin{lemma}{section}
\numberwithin{whatsoever}{section}
\numberwithin{equation}{section}

\begin{document}
\title{Realization of HCMU metrics in  3-dimensional space forms as Weingarten surfaces}
\author{Zhiqiang Wei}
\address{School of Mathematical and Statistics \\
Henan University \\
Kaifeng 475004 P.R. China}
\email{weizhiqiang15@mails.ucas.edu.cn}
\author{Yingyi Wu}
\address{School of Mathematical Sciences \\
University of Chinese Academy of Sciences \\
Beijing 100049 P.R. China  }
\email{wuyy@ucas.ac.cn}

\begin{abstract}  In  \cite{WeiWu22-3} (J. Geom. Anal. 32, 199(2022)), we classified HCMU surfaces in 3-dimensional Euclidean space as Weingarten surfaces by holomorphic functions. In this paper, using a totally different method  from  \cite{WeiWu22-3},  we  will  classify HCMU surfaces in 3-dimensional  space  forms as Weingarten surfaces. Moreover we will give a criteria of a Weingarten surface which is also an  HCMU surface.

\smallskip
2010 MSC 53A50, 32Q15.

\end{abstract}

\maketitle

%\tableofcontents

\section{Introduction}
\label{section: introduction}

~~~~~In a fixed K\"{a}hler class on a compact K$\ddot{a}$hler manifold $M$, an extremal K\"{a}hler metric which was introduced by Calabi \cite{Ca} is the critical
point of the following Calabi energy functional
$$
\mathcal{C}(g)=\int_{M} R^2 dg,
$$
where $R$ is the scalar curvature of the metric $g$ in the given K$\ddot{a}$hler class. The Euler-Lagrange equations of $\mathcal{C}(g)$ are $R_{,\alpha
\beta}=0$ for all indices $\alpha, \beta$, where $R_{,\alpha\beta}$ is the second-order $(0,2)$
covariant derivative of $R$. When $M$ is a compact Riemann surface,
 Calabi in \cite{Ca} proved that an extremal K$\ddot{a}$hler metric is a CSC (constant scalar curvature) metric. \par

A natural question is whether or not an extremal K$\ddot{a}$hler metric with singularities on a compact Riemann surface is still a CSC metric. In \cite{Ch1}, Chen  first gave an example of a non-CSC extremal K\"{a}hler metric with singularities. Now we often call a non-CSC extremal K$\ddot{a}$hler metric with finite singularities on a compact Riemann surface an HCMU (the Hessian of the Curvature of the Metric is Umbilical) metric.\par

Since \cite{Ch1}, HCMU metrics  has been widely studied. In \cite{Ch2}, Chen gave an obstruction theorem of the existence of HCMU metrics. In \cite{LZ}, Lin and Zhu proved that the complex gradient of Gauss curvature $K$ (i.e., $\nabla K=\sqrt{-1}e^{-2\varphi }K_{\overline{z}}\frac{\partial}{\partial z}$ when the metric is $g=e^{2\varphi}|dz|^{2}$) is a meromorphic vector field on the underlying Rieamnn surface and gave an explicit construction of a kind of HCMU metrics on $S^{2}$. In \cite{Ch2} and \cite{LZ}, Chen and Lin-Zhu independently proved that the Gauss curvature of an HCMU can be continuously extended to the whole underlying compact Riemann surface. In \cite{CCW}, Chen, Chen and the second author proved that any HCMU metric is locally isometric to an HCMU metric on $S^{2}$ with only two singularities (they called this kind of HCMU metrics football). In \cite{CW2}, Chen and the second author reduced the existence of an HCMU metric with conical singularities to the existence of a meromorphic 1-form on the underlying Riemann surface. Later, in \cite{Xb}, Chen, Xu and the second author reduced the existence of an HCMU metric with both conical singularities and cusp singularities to the existence of a meromorphic 1-form on the underlying Riemann surface. It is interesting that on any compact Riemann surface there always exists an HCMU metric without considering the positions of singularities.  However, if considering the positions of the singularities, the existence of  HCMU metrics is very complication. More  results of the existence of HCMU metrics, one can refer to\cite{Ch2},\cite{CW1},\cite{WeiWu18} and the references cited in.\par

Recently, Peng and the second author \cite{PW} proposed the question: whether or not an HCMU metirc can be isometrically imbedded into some ``good" spaces, for example, 3-dimensional space forms or $\mathbb{C}P^{N}$. For $\mathbb{R}^{3}$, they  \cite{PW} got a one-parameter family of immersions into $\mathbb{R}^{3}$, each of whom is a Weingarten surface. In \cite{WeiWu22-1}, we proved that any  HCMU metric locally can be isometrically imbedded into 3-dimensional space forms as Weingarten surfaces. As an application, we proved that any HCMU metric locally can be isometrically imbedded into $\mathbb{C}P^{3}$. In \cite{WeiWu22-3}, we classified the isometrically imbedded of an HCMU metric into $\mathbb{R}^{3}$ as Weingarten surfaces by holomorphic functions.  In \cite{WeiWu22-4}, we proved that any  HCMU metric locally can not be  minimal immersed into 3-dimensional space forms by moving frames. \par

In this paper, we will systematic study the realization of an HCMU metric into 3-dimensional space forms as Weingarten surfaces by Hopf differential.\par

Let $g$ be an HCMU metric on a compact Riemann surface $M$ with  the character 1-form $\omega$. Suppose the Gauss curvature of $g$ is $K$ with the maximum $K_{1}$ and the minimum $K_{2}$. Denote $-\frac{1}{3}(K-K_{1})(K-K_{2})(K+K_{1}+K_{2})$ by $p$. Then $g=4p \omega\overline{\omega}$ (see \cite{CW2}). Denote $M^{*}=M\setminus\{ \text{zeros and poles of} ~\omega\}$. Let $\mathbb{Q}^{3}_{c}$ be the space form of constant sectional curvature $c$.  Our main results can be stated as  follows.\par

\begin{theorem}\label{mainth-1}
For any point $P\in M^{*}$, there exist a neighborhood $U\subseteq M^{*}$ of $P$ and an isometric immersion $F:U\rightarrow Q^{3}_{c}$ which is a Weingarten surface without umbilical points if and only if there exists a real constant $A$ such that the mean curvature $H$ of $F$ satisfying
\begin{equation*}
H'=\frac{p-2p'(H^{2}-K+c)}{2[pH\pm\sqrt{p^{2}(H^{2}-K+c)-A^{2}}~]},
\end{equation*}
where `` $ ' $ " denotes ``$\frac{d}{dK}$".
\end{theorem}

\begin{theorem}\label{mainth-2}
Let $z=x+\sqrt{-1}y$ be a complex coordinate on $U\subseteq M^{*}$ such that the HCMU metric of
 $r:U\rightarrow \mathbb{Q}^{3}_{c}$ can be written as
$$g=-\frac{4}{3}(K-K_{1})(K-K_{2})(K+K_{1}+K_{2})|dz|^{2}.$$
Suppose the second fundamental form of $r$ is
$$B=h_{11}dx^{2}+2h_{12}dxdy+h_{22}dy^{2}.$$
Then $r$ is a Weingarten surface if and only if $h_{12}$ is a constant.
\end{theorem}

\section{Preliminaries }
\begin{Definition}[\cite{T}]\label{conedef}
 Let $M$ be a Riemann surface, $P\in M$. A conformal metric $g$ on $M$ is said to have a conical singularity at
 $P$ with the singular angle $2 \pi \alpha~(\alpha>0,\alpha\neq 1)$ if in a neighborhood of $P$
 \begin{equation*}
  g=e^{2\varphi}|dz|^2,
 \end{equation*}
 where $z$ is a local complex coordinate defined in the neighborhood of $P$ with $z(P)=0$ and
 \begin{equation*}
  \varphi-(\alpha-1)\ln|z|
 \end{equation*}
 is continuous at $0$.
\end{Definition}
\begin{Definition}[\cite{Xb}]\label{cuspdef}
 Let $M$ be a Riemann surface, $P\in M$. A conformal metric $g$ on $M$ is said to have a cusp singularity at
 $P$ if in a neighborhood of $P$
 \begin{equation*}
  g=e^{2\varphi}|dz|^2,
 \end{equation*}
 where $z$ is a local complex coordinate defined in the neighborhood of $P$ with $z(P)=0$ and
 \begin{equation*}
 \lim_{z\rightarrow 0} \frac{\varphi+\ln|z|}{\ln|z|}=0.
 \end{equation*}
\end{Definition}
\begin{Definition}[\cite{Ch2}]
 Let $M$ be a compact Riemann surface and $P_1,\cdots,P_N$ be $N$ points on $M$.
 Denote $M\backslash \{P_1,\ldots,P_N\}$ by $M^*$. Let $g$ be a conformal metric on $M^*$.
 If the Gauss curvature $K$ of $g$ satisfies
\begin{equation}\label{HCMUequ}
  K_{,zz}=0
\end{equation}
and $K$ is not a constant, we call $g$ an HCMU metric on $M$.
\end{Definition}
  In this paper we always consider HCMU metrics with finite area and
  finite Calabi energy, that is,
\begin{equation*}
 \int_{M^*}dg<+\infty, ~~ \int_{M^*}
K^2 dg <+\infty.
\end{equation*}
\par
From \cite{Ch0}, \cite{LZ}, \cite{WZ},  we know that each singularity of an HCMU metric  is  conical singularity or cusp if it has finite area and finite Calabi energy.\par
We now list some results of HCMU metrics, which will be used in this paper. For more results one can refer to \cite{CW2},\cite{Xb} and the references cited in.\par
First the equation (\ref{HCMUequ}) is equivalent to
\begin{equation*}
 \nabla K =\sqrt{-1}e^{-2\varphi}K_{\bar{z}}\frac{\partial}{\partial z},
\end{equation*}
which is a holomorphic vector field on $M^*$. In \cite{LZ}, Lin and Zhu proved that $\nabla K$ is actually a meromorphic vector field on $M$.
In \cite{CW2}, Chen and the second author defined the dual 1-form of $\nabla K$ by $\omega(\nabla K)=\frac{\sqrt{-1}}{4}$.
They call $\omega$ the character 1-form of the metric. In \cite{Ch2},\cite{LZ}, the authors proved that the curvature
$K$ can be continuously extended to $M$ and there are finite smooth extremal points of $K$ on
$M^*$. In \cite{CCW},\cite{Xb}, the authors proved the following facts: each smooth extremal
point of $K$ is either the maximum point of $K$ or the minimum point of $K$, and if we denote the maximum
of $K$ by $K_1$ and the minimum of $K$ by $K_2$ then if all the singularities of $g$ are conical singularities,
$$
 K_1>0,~K_1>K_2>-(K_1+K_2);
$$
if there exist cusps in the singularities,
$$K_{1}>0,~K_{2}=-\frac{1}{2}K_{1}.$$
 Denote $M^* \setminus \{\text{smooth extremal
points of}~K \}$ by $M'$. Then on $M'$
\begin{equation}\label{sys0}
\begin{cases}
 \cfrac{dK}{-\frac{1}{3}(K-K_1)(K-K_2)(K+K_1+K_2)}=\omega+\bar{\omega}, \\
g=-\frac{4}{3}(K-K_1)(K-K_2)(K+K_1+K_2)\omega \bar{\omega}.\\
\end{cases}
\end{equation}
By (\ref{sys0}), some properties of $\omega$ are got in \cite{CW2}:
\begin{itemize}
 \item $\omega$ only has simple poles,
 \item at each pole, the residue of $\omega$ is a non-zero real number,
  \item $\omega+\bar{\omega}$ is exact on $M \setminus \{poles~of~\omega\}$.
 \end{itemize}
Conversely, if a meromorphic 1-form $\omega$ on $M$ which satisfies the properties above, then we pick two real numbers $K_{1},K_{2}$ such that $K_{1}>0,K_1>K_2>-(K_1+K_2)$ or $K_{1}>0,~K_{2}=-\frac{1}{2}K_{1},$  and consider the following equation on $M\setminus\{\text{poles of}~\omega\}$
\begin{equation}\label{sys1}
\begin{cases}
 \cfrac{dK}{-\frac{1}{3}(K-K_1)(K-K_2)(K+K_1+K_2)}=\omega+\bar{\omega}, \\
K(P_{0})=K_{0},\\
\end{cases}
\end{equation}
where $P_{0}\in M\setminus\{\text{poles of}~\omega\}$ and $K_{2}< K_{0}< K_{1}$. We get that (\ref{sys1}) has a unique solution $K$ on $M\setminus\{\text{poles of}~\omega\}$ and $K$ can be continuously extended to $M$. Furthermore, we define a metric $g$ on $M\setminus\{\text{poles of}~\omega\}$ by
$$g=-\frac{4}{3}(K-K_{1})(K-K_{2})(K+K_{1}+K_{2})\omega\overline{\omega},$$
 where $K$ is the solution of (\ref{sys1}). Then it can be proved that $g$ is an HCMU metric, $K$ is the Gauss curvature of $g$ and $\omega$ is the character 1-form of $g$.\par
It is interesting that on any compact Riemann surface there always exists a meromorphic 1-form  satisfying the properties (see \cite{CW2}). So there always exists an HCMU metric on a compact Riemann surface.

\section{Some known results}
In this section, we will recall some results of isometric immersion of HCMU metrics into 3-dimensional space forms. For more results, one may refer to \cite{PW} and \cite{ WeiWu22-1,WeiWu22-3}. \par
Suppose $g$ is an HCMU metric on compact Riemann surface $M$ with  the character 1-form $\omega$. Let $K$ be the Gauss curvature of $g$ with the maximum $K_{1}$ and the minimum $K_{2}$. Denote $M^{*}=M\setminus\{ \text{zeros and poles of} ~\omega\}$. Let $\mathbb{Q}^{3}_{c}$ be the 3-dimensional space form of constant sectional curvature $c$. In \cite{WeiWu22-1}, we  proved the following theorem.

\begin{Theorem}[\cite{WeiWu22-1}]\label{WeiWu22-1}
For any point $P\in M^{*}$, there exist a  neighborhood $U\subseteq M^{*}$ of $P$ and an isometric immersion $F:U\rightarrow \mathbb{Q}^{3}_{c}$ which is a Weingarten surface.
\end{Theorem}
When $\mathbb{Q}^{3}_{c}=\mathbb{R}^{3}$, Peng and the second author using moving frames proved the following theorem.

\begin{Theorem}[\cite{PW}]
Denote $\sqrt{-\frac{4}{3}(K-K_{1})(K-K_{2})(K+K_{1}+K_{2})}$ by $\mu(K)$ and $\frac{1}{6}(K_{1}^{2}+K_{1}K_{2}+K_{2}^{2})$ by $\delta$. Then at any point $P$ of $M^{*}$ there is  an open neighborhood $U$ of $P$ such that exists a one-parameter family of isometric immersions of $g$ from $U$
 to $\mathbb{R}^{3}$ such the mean curvature of each immersion has the following expression:
 $$H=\pm \frac{1}{\mu(K)}(\frac{-K\mu^{2}(K)}{4\sqrt{-\frac{K^{4}}{4}+\delta K^{2}+s}}-\sqrt{-\frac{K^{4}}{4}+\delta K^{2}}+s),$$
 where $s$ is a constant with $s\geq \frac{1}{12}K_{1}^{2}(K_{1}^{2}-2K_{1}K_{2}-2K_{2}^{2})$.
 \end{Theorem}

We  \cite{WeiWu22-3} classified HCMU surfaces as Weingarten surfaces in $\mathbb{R}^{3}$ by holomorphic functions.

\begin{Theorem}[\cite{WeiWu22-3}]\label{WeiWu22-2}
On $M^{*}$, if $\mathcal{S}$ is a local isometric immersion into $\mathbb{R}^{3}$, that is, $\mathcal{S}$ is an HCMU surface and $\mathcal{S}$ is a Weingarten surface without umbilical points, then the mean curvature $H$ of $\mathcal{S}$ satisfies
\begin{equation}\label{E-T-1}
H'=\frac{2p'(K)(H^{2}-K)-p(K)}{\sqrt{4p^{2}(K)(H^{2}-K)-\varepsilon}-2p(K)H},
\end{equation}
or
\begin{equation}\label{E-T-2}
H'=-\frac{2p'(K)(H^{2}-K)-p(K)}{\sqrt{4p^{2}(K)(H^{2}-K)-\varepsilon}+2p(K)H},
\end{equation}
where ``~$'$" means ``~$\frac{d}{dK}$", $p(K)=-\frac{1}{3}(K-K_{1})(K-K_{2})(K+K_{1}+K_{2})$ and $\varepsilon$ is a constant with $\varepsilon\geq0$ and $4p^{2}(K)(H^{2}-K)-\varepsilon>0$. Conversely, if $H$ is a solution of (\ref{E-T-1}) or (\ref{E-T-2}), then there exists a local HCMU surface $\mathcal{S}$ which is a Weingarten surface in $\mathbb{R}^{3}$ without umbilical points such that the mean curvature of $\mathcal{S}$ is $H$.
\end{Theorem}

\section{Proofs of Theorem  \ref{mainth-1} and \ref{mainth-2}}
Suppose $g$ is an HCMU metric on compact Riemann surface $M$ with  the character 1-form $\omega$. Let $K$ be the Gauss curvature of $g$ with the maximum $K_{1}$ and the minimum $K_{2}$. Denote $M^{*}=M\setminus\{ \text{zeros and poles of} ~\omega\}$. Let $\mathbb{Q}^{3}_{c}$ be the 3-dimensional space form of constant sectional curvature $c$.
Suppose $F:U\rightarrow \mathbb{Q}^{3}_{c}$ is an isometric immersion without umbilical point, where $P\in M^{*}$ and $ U\subseteq M^{*}$  is an open neighborhood of $P$. Since $\omega$ is a non-vanishing holomorphic 1-form on $M^{*}$, locally we can write  $\omega$ as $dz$, where $z=x+\sqrt{-1}y$ is a local complex coordinate defined in a neighborhood of $P$. Moreover
$$g=e^{u}|dz|^{2}=-\frac{4}{3}(K-K_{1})(K-K_{2})(K+K_{1}+K_{2})|dz|^{2},$$
where $e^{u}= -\frac{4}{3}(K-K_{1})(K-K_{2})(K+K_{1}+K_{2})$.\par
 Denote $-\frac{1}{3}(K-K_{1})(K-K_{2})(K+K_{1}+K_{2})$ by $p(K)$. Then by equation (\ref{sys1})
\begin{equation*}
\frac{dK}{p(K)}=2dx,
\end{equation*}
so $K$  depends only  on $x$. In fact,
\begin{equation}\label{K-1}
x=\frac{1}{2}\ln[(K_{1}-K)^{\sigma}(K-K_{2})^{\beta}(K+K_{1}+K_{2})^{\gamma}]+c_{1},
\end{equation}
where $\sigma= \frac{-3}{(K_{1}-K_{2})(K_{2}+2K_{1})}, \beta= \frac{-3}{(K_{2}-K_{1})(2K_{2}+K_{1})}, \gamma=  \frac{-3}{(K_{2}+2K_{1})(2K_{2}+K_{1})}$, $c_{1}$ is a real constant.\par
Suppose that the second fundamental form of $F:U\rightarrow \mathbb{Q}^{3}_{c}$ is
\begin{equation*}
B=h_{11}dx^{2}+2h_{12}dxdy+h_{22}dy^{2}.
\end{equation*}
Then the mean curvature $H$ of $F$ is
\begin{equation}\label{Mean-C}
H=\frac{1}{2e^{u}}(h_{11}+h_{22}).
\end{equation}

Let
\begin{equation}\label{Hopf-D}
Q=\frac{1}{4}(h_{11}-h_{22}-2\sqrt{-1}h_{12}).
\end{equation}
Then $Qdz\otimes dz$ is a globally defined quadratic form which is called the Hopf differential of the surface. Furthermore, the Gauss equation can be written as (see \cite{CL})
\begin{equation*}
4|Q|^{2}=e^{2u}(H^{2}-K+c).
\end{equation*}

Suppose that
\begin{equation*}
Q=\frac{1}{2}e^{u}\sqrt{H^{2}-K+c}(\cos\theta+\sqrt{-1}\sin\theta).
\end{equation*}
Since $F$ has no umbilical point, $\theta$ is a smooth function of $x,y$.\par
By equations (\ref{Mean-C}) and (\ref{Hopf-D}), we obtain
\begin{equation*}
\begin{cases}
h_{11}=e^{u}H+e^{u}\sqrt{H^{2}-K+c}\cos\theta,\\
h_{12}=-e^{u}\sqrt{H^{2}-K+c}\sin\theta,\\
h_{22}=e^{u}H-e^{u}\sqrt{H^{2}-K+c}\cos\theta.
\end{cases}
\end{equation*}

If $F(U)$ is a Weingarten surface without umbilical point, i.e., $dK\wedge dH=0$, we obtain $H_{y}=0$, that is, $H$ also depends only on $x$.

The Codazzi equation (see \cite{CL})
\begin{equation*}\label{Codazzi-0}
Q_{\overline{z}}=\frac{1}{2}e^{u}H_{z},
\end{equation*}
can be written as
\begin{equation}\label{Codazzi-1}
\begin{cases}
H_{x}\sin\theta=-\sqrt{H^{2}-K+c}\theta_{x},\\
e^{u}H_{x}\cos\theta=-e^{u}\sqrt{H^{2}-K+c}\theta_{y}+(e^{u}\sqrt{H^{2}-K+c})_{x},
\end{cases}
\end{equation}

\begin{Proposition}
$H$ is not a constant.
\end{Proposition}

\begin{proof}
If $H$ is a constant, since $K$ is not a constant, then $H^{2}-K+c$  is not identically vanishing. Thus the Codazzi equation  (\ref{Codazzi-1}) converts to
\begin{equation*}
\begin{cases}
\theta_{x}=0,\\
\theta_{y}=(\ln e^{u}\sqrt{H^{2}-K+c})_{x}.
\end{cases}
\end{equation*}
So there exists a constant $A\neq0$ such that
\begin{equation*}
\begin{cases}
\theta_{x}=0,\\
\theta_{y}=A,\\
(\ln e^{u}\sqrt{H^{2}-K+c})_{x}=A.
\end{cases}
\end{equation*}
Since $\frac{dK}{dx}=2p$,

\begin{equation*}
(\ln e^{u}\sqrt{H^{2}-K+c})_{K}=\frac{A}{2p}.
\end{equation*}
By
\begin{equation*}
\frac{d}{dK}\ln [(K_{1}-K)^{\sigma}(K-K_{2})^{\beta}(K+K_{1}+K_{2})^{\gamma}]=\frac{1}{p(K)},
\end{equation*}
where $\sigma= \frac{-3}{(K_{1}-K_{2})(K_{2}+2K_{1})} ,\beta= \frac{-3}{(K_{2}-K_{1})(2K_{2}+K_{1})} ,\gamma=  \frac{-3}{(K_{2}+2K_{1})(2K_{2}+K_{1})}$, we obtain
\begin{equation*}
\begin{cases}
\theta=Ay+A_{1},\\
\ln e^{u}\sqrt{H^{2}-K+c}=\frac{A}{2}\ln [(K_{1}-K)^{\sigma}(K-K_{2})^{\beta}(K+K_{1}+K_{2})^{\gamma}]+A_{2},
\end{cases}
\end{equation*}
where $A_{1},A_{2}$ are constants.
Thus
\begin{equation*}
e^{u}\sqrt{H^{2}-K+c}= e^{A_{2}}[(K_{1}-K)^{\sigma}(K-K_{2})^{\beta}(K+K_{1}+K_{2})^{\gamma}]^{\frac{A}{2}}.
\end{equation*}

Since $e^{u}=4p(K)$,
$$H^{2}=\frac{e^{2A_{2}}}{16p^{2}}[(K_{1}-K)^{\sigma}(K-K_{2})^{\beta}(K+K_{1}+K_{2})^{\gamma}]^{A}+K-c.$$

This is a contradiction!

\end{proof}

\begin{Remark}
In \cite{WeiWu22-4}, using moving frames, we have proved this property.
\end{Remark}

\begin{Proposition}
$H^{2}-K+c$ is not identically vanishing.
\end{Proposition}

\begin{proof}
Suppose $H^{2}-K+c\equiv0$. The Codazzi equation (\ref{Codazzi-1}) can be written as
\begin{equation*}
\begin{cases}
H_{x}\sin\theta=0,\\
H_{x}\cos\theta=0.
\end{cases}
\end{equation*}
Then $H$ is a constant. It is a contradiction.
\end{proof}

Now the Codazzi equation (\ref{Codazzi-1}) can be written as
\begin{equation}\label{Codazzi-2}
\begin{cases}
\theta_{x}=-\frac{H_{x}}{\sqrt{H^{2}-K+c}}\sin\theta,\\
\theta_{y}=(\ln e^{u}\sqrt{H^{2}-K+c})_{x}-\frac{H_{x}}{\sqrt{H^{2}-K+c}}\cos\theta.
\end{cases}
\end{equation}

\begin{Proposition}
$\theta_{y}=0$, i.e., $\theta$ is a function  depending only on $x$.
\end{Proposition}

\begin{proof}
By the Codazzi equation (\ref{Codazzi-2}),
\begin{equation}\label{Codazzi-3}
\begin{cases}
\theta_{xy}=-\frac{H_{x}}{\sqrt{H^{2}-K+c}}\cos\theta \theta_{y},\\
\theta_{yx}=(\ln e^{u}\sqrt{H^{2}-K+c})_{xx}-(\frac{H_{x}}{\sqrt{H^{2}-K+c}})_{x}\cos\theta+\frac{H_{x}}{\sqrt{H^{2}-K+c}}\sin\theta \theta_{x}.
\end{cases}
\end{equation}

Since $\theta_{xy}=\theta_{yx}$,
\begin{equation*}
\frac{H_{x}}{\sqrt{H^{2}-K+c}}(\sin\theta \theta_{x}+\cos\theta \theta_{y})=(\frac{H_{x}}{\sqrt{H^{2}-K+c}})_{x}\cos\theta-(\ln e^{u}\sqrt{H^{2}-K+c})_{xx}.
\end{equation*}
Also by the Codazzi equation (\ref{Codazzi-2}),
\begin{equation*}
\sin\theta \theta_{x}+\cos\theta \theta_{y}=(\ln e^{u}\sqrt{H^{2}-K+c})_{x}\cos\theta- \frac{H_{x}}{\sqrt{H^{2}-K+c}}.
\end{equation*}
Thus we obtain
\begin{equation*}
\begin{aligned}
&[(\frac{H_{x}}{\sqrt{H^{2}-K+c}})_{x}-\frac{H_{x}}{\sqrt{H^{2}-K+c}} \times (\ln e^{u}\sqrt{H^{2}-K+c})_{x}]\cos\theta\\
&~~~~~~~~~~~~=(\ln e^{u}\sqrt{H^{2}-K+c})_{xx}-\frac{H^{2}_{x}}{H^{2}-K+c}.
\end{aligned}
\end{equation*}
Since $H,K$ are functions depending only on $x$, the following two cases may happen:\par
\textbf{Case A:}
\begin{equation}\label{case-2}
\begin{cases}
(\frac{H_{x}}{\sqrt{H^{2}-K+c}})_{x}-\frac{H_{x}}{\sqrt{H^{2}-K+c}} \times (\ln e^{u}\sqrt{H^{2}-K+c})_{x}=0,\\
(\ln e^{u}\sqrt{H^{2}-K+c})_{xx}-\frac{H^{2}_{x}}{H^{2}-K+c}=0.
\end{cases}
\end{equation}
\par
\textbf{Case B:} $\theta$ is a function depending only on $x$, i.e., $\theta_{y}=0$.\par

Suppose \textbf{Case A} holds. By the first equation in (\ref{case-2}), there exists a constant $A\neq0$ such that
\begin{equation}\label{H-1}
H_{x}= Ae^{u}(H^{2}-K+c).
\end{equation}

Then by the second equation in (\ref{case-2}),
\begin{equation}\label{case-2-2}
 (\ln e^{u}\sqrt{H^{2}-K+c})_{xx}=A^{2}e^{2u}(H^{2}-K+c).
\end{equation}

Set
\begin{equation*}
 e^{u}\sqrt{H^{2}-K+c}=e^{\xi(x)}.
\end{equation*}
then equation (\ref{case-2-2}) can be written as
\begin{equation}\label{ODE-1}
\xi''=A^{2}e^{2\xi},
\end{equation}
where `` $ ' $ " denotes`` $ \frac{d}{dx} $ ".\par
Multiplying $2\xi'$ in both sides of (\ref{ODE-1}), we obtain
\begin{equation*}
[(\xi')^{2}]'=2\xi'\xi''=2A^{2}e^{2\xi}\xi'=A^{2}(e^{2\xi})',
\end{equation*}
so there exists a constant $B$ such that
\begin{equation*}
(\xi')^{2}=A^{2}(e^{2\xi}+B),
\end{equation*}
i.e.,
\begin{equation*}
\xi'=\pm A\sqrt{e^{2\xi}+B }.
\end{equation*}
Then
\begin{equation*}
\frac{de^{2\xi}}{e^{2\xi}\sqrt{e^{2\xi}+B}}=\pm 2Adx.
\end{equation*}

Set $\sqrt{e^{2\xi}+B}=t,$ then $e^{2\xi}=t^{2}-B$. Then
\begin{equation*}
\frac{dt}{t^{2}-B}=\pm Adx.
\end{equation*}

(a) If $B=0$,
\begin{equation*}
t=\frac{1}{\pm Ax+D_{1}},
\end{equation*}
where $D_{1}$ is a constant. \par
By $e^{2u}(H^{2}-K+c)=e^{2\xi}=t^{2}$,
\begin{equation*}
H=\pm\sqrt{\frac{1}{16p^{2}}\cdot\frac{1}{(\pm Ax+D_{1})^{2}}+K-c}.
\end{equation*}

However, by (\ref{H-1}),
\begin{equation*}
H_{x}= \frac{Ae^{2\xi}}{e^{u}}=\frac{A t^{2}}{e^{u}}=\frac{A}{4p(\pm Ax+D_{1})^{2}},
\end{equation*}

Since (\ref{K-1}), we obtain a contradiction.

(b) If $B>0$
\begin{equation*}
t=\frac{\sqrt{B}(e^{\pm 2A\sqrt{B}x+D_{1}}+1)}{1-e^{\pm 2A\sqrt{B}x+D_{1}}},
\end{equation*}
where $D_{1}$ is a constant. Thus,
\begin{equation*}
H=\pm\sqrt{\frac{Be^{\pm 2A\sqrt{B}x+D_{1}}}{4p^{2}(1-e^{\pm 2A\sqrt{B}x+D_{1}})^{2}}+K-c}.
\end{equation*}

(c) If $B<0$,
\begin{equation*}
t=\sqrt{-B}\tan [\sqrt{-B}(\pm Ax+D_{1})],
\end{equation*}
where $D_{1}$ is a constant. Thus,

\begin{equation*}
H=\pm \sqrt{\frac{-B(\tan^{2}[\sqrt{-B}(\pm Ax+D_{1})]+1)}{16p^{2}}+K-c}.
\end{equation*}

Similar as case (a), for cases (b)(c), we obtain $H_{x}\neq \pm Ae^{u}(H^{2}-K+c)$, so $\theta$ is a function depending only on $x$, i.e., $\theta_{y}=0$.
\end{proof}

\begin{Remark}
By $\frac{dK}{dx}=2p$, equation (\ref{case-2}) converts to
\begin{equation*}
\begin{cases}
H'=\pm 2A(H^{2}-K+c),\\
(2p''+4p'HA)(H^{2}-K+c)^{2}+(4AH-p')(H^{2}-K+c)-p=0,
\end{cases}
\end{equation*}
where `` $ ' $ " denotes ``$\frac{d}{dx}$". So the above equation has no solution.
\end{Remark}

\begin{Proposition}\label{Pro-1}
There exists a constant $A$ such that
\begin{equation}\label{Thm-e-1}
H'=\frac{p-2p'(H^{2}-K+c)}{2[pH\pm\sqrt{p^{2}(H^{2}-K+c)-A^{2}}]},
\end{equation}
where `` $ ' $ " denotes ``$\frac{d}{dK}$".
\end{Proposition}
\begin{proof}
Since $\theta$ is a function depending only on $x$, the Codazzi equation (\ref{Codazzi-2}) converts to
\begin{equation}\label{Codazzi-x}
\begin{cases}
\frac{H_{x}}{\sqrt{H^{2}-K+c}}\sin\theta=-\theta_{x},\\
\frac{H_{x}}{\sqrt{H^{2}-K+c}}\cos\theta=(\ln e^{u}\sqrt{H^{2}-K+c})_{x}.
\end{cases}
\end{equation}

(1)\textbf{ If $\theta\equiv0$ },  the Codazzi equation (\ref{Codazzi-x}) converts to
$$\frac{H_{x}}{\sqrt{H^{2}-K+c}}=(\ln e^{u}\sqrt{H^{2}-K+c})_{x},$$
i.e.,
$$\frac{H'}{\sqrt{H^{2}-K+c}}=(\ln e^{u}\sqrt{H^{2}-K+c})'.$$
where `` $ ' $ " denotes ``$\frac{d}{dK}$".
That is
$$H'=\frac{p-2p'(H^{2}-K+c)}{2p[H-\sqrt{H^{2}-K+c}]}.$$

(2) \textbf{If $\theta\equiv\pi$},  the Codazzi equation (\ref{Codazzi-x}) converts to
$$\frac{H_{x}}{\sqrt{H^{2}-K+c}}=-(\ln e^{u}\sqrt{H^{2}-K+c})_{x}.$$
Similar as the case (1),
$$H'=\frac{p-2p'(H^{2}-K+c)}{2p[H+\sqrt{H^{2}-K+c}]},$$
where `` $ ' $ " denotes ``$\frac{d}{dK}$"

(3) \textbf{If $\sin\theta\neq0$}, by the Codazzi equation (\ref{Codazzi-x}),
$$-\frac{\cos\theta \theta_{x}}{\sin\theta}=-\frac{d\ln|\sin\theta|}{dx}=(\ln e^{u}\sqrt{H^{2}-K+c})_{x},$$
i.e.,
$$[\ln |\sin\theta|e^{u}\sqrt{H^{2}-K+c} ]_{x}=0,$$
so there exists a constant $A\neq0$ such that (Noticing $e^{u}=4p$)
$$\sin\theta= \frac{A}{p\sqrt{H^{2}-K+c}}\in [-1,1].$$
Then
$$\cos\theta=\pm\frac{\sqrt{p^{2}(H^{2}-K+c)-A^{2}}}{p\sqrt{H^{2}-K+c}},$$
$$\cos\theta\theta'=-\frac{A[2p'(H^{2}-K+c)+p(2HH'-1)]}{2p^{2}(H^{2}-K+c)\sqrt{H^{2}-K+c}},$$
where `` $ ' $ " denotes ``$\frac{d}{dK}$.\par
Thus, if
\begin{equation*}
\begin{cases}
\sin\theta= \frac{A}{p\sqrt{H^{2}-K+c}},\\
\cos\theta=\frac{\sqrt{p^{2}(H^{2}-K+c)-A^{2}}}{p\sqrt{H^{2}-K+c}},\\
\theta'=-\frac{A[2p'(H^{2}-K+c)+p(2HH'-1)]}{2p(H^{2}-K+c)\sqrt{p^{2}(H^{2}-K+c)-A^{2}}}
\end{cases}
\end{equation*}
the Codazzi equation (\ref{Codazzi-x}) converts to
$$H'=\frac{p-2p'(H^{2}-K+c)}{2[pH-\sqrt{p^{2}(H^{2}-K+c)-A^{2}}]},$$
where `` $ ' $ " denotes ``$\frac{d}{dK}$".

If
\begin{equation*}
\begin{cases}
\sin\theta= \frac{A}{p\sqrt{H^{2}-K+c}},\\
\cos\theta=-\frac{\sqrt{p^{2}(H^{2}-K+c)-A^{2}}}{p\sqrt{H^{2}-K+c}},\\
\theta'=\frac{A[2p'(H^{2}-K+c)+p(2HH'-1)]}{2p(H^{2}-K+c)\sqrt{p^{2}(H^{2}-K+c)-A^{2}}}
\end{cases}
\end{equation*}
the Codazzi equation (\ref{Codazzi-x}) converts to
$$H'=\frac{p-2p'(H^{2}-K+c)}{2[pH+\sqrt{p^{2}(H^{2}-K+c)-A^{2}}~]},$$
where `` $ ' $ " denotes ``$\frac{d}{dK}$".

To sum up, there exists a constant $A$ such that
$$H'=\frac{p-2p'(H^{2}-K+c)}{2[pH\pm\sqrt{p^{2}(H^{2}-K+c)-A^{2}}~]},$$
where `` $ ' $ " denotes ``$\frac{d}{dK}$".
\end{proof}
By the fundamental theorem of hypersurfaces (see \cite{MR}) in space forms and \textbf{Proposition \ref{Pro-1}}, we finish the proof of \textbf{Theorem \ref{mainth-1}}.

\begin{Remark}
If $A=0$, by direct calculation,
\begin{equation*}
\begin{aligned}
H=&\mp\frac{1}{2\sqrt{p}}[\frac{(c-K)p(K)}{\sqrt{-\frac{1}{4}K^{4}+\frac{c}{3}K^{3}+\delta K^{2}-2c\delta K+s }}-\sqrt{-\frac{1}{4}K^{4}+\frac{c}{3}K^{3}+\delta K^{2}-2c\delta K+s }~],
\end{aligned}
\end{equation*}
is the solution of equation (\ref{Thm-e-1}), where $\delta=\frac{1}{6}(K_{1}^{2}+K_{1}K_{2}+K_{2}^{2})$, and $s$ is a constant.
\end{Remark}

Now we give more results of realization of HCMU metrics in $\mathbb{Q}^{3}_{c}$.
\begin{Corollary}
$h_{12}$ is a constant.
\end{Corollary}
\begin{proof}
By the proof of \ref{Pro-1}, there exists a constant $A$ such that
$$e^{u}\sqrt{H^{2}-K+c}\sin\theta =A,$$
so
$$h_{12}=-e^{u}\sqrt{H^{2}-K+c}\sin\theta=-A.$$
\end{proof}

\begin{Theorem}\label{Thm-1}
Let $r:U \rightarrow \mathbb{Q}^{3}_{c}$ be an HCMU surface with metric
$$g=-\frac{4}{3}(K-K_{1})(K-K_{2})(K+K_{1}+K_{2})|dz|^{2}.$$
Suppose the second fundamental form of $r$ is $$B=h_{11}dx^{2}+2h_{12}dxdy+h_{22}dy^{2}.$$
If $h_{12}$ is a constant, then $r$ is a Weingarten surface.
\end{Theorem}
\begin{proof}
Since $h_{12}$ is a constant, the Codazzi equation of $r$ is
\begin{equation*}
\begin{cases}
h_{11y}=0,\\
(h_{11}+h_{22})u_{x}=2h_{22x}.
\end{cases}
\end{equation*}
Then we can suppose $h_{11}=f(x)$.  By the Gauss equation
$$e^{2u}(K-c)=h_{11}h_{22}-h_{12}^{2},$$
$f\neq0$, and
$$h_{22}=\frac{e^{2u}(K-c)+h_{12}^{2}}{f(x)}.$$
So $h_{22}$ is a function depending only on $x$.  Moreover, $H=\frac{h_{11}+h_{22}}{2e^{u}}$ is a function depending only on $x$. Then $r:\Sigma\rightarrow \mathbb{Q}^{3}_{c}$ is a Weingarten surface.
\end{proof}

By \textbf{Proposition \ref{Pro-1}} and \textbf{Theorem \ref{Thm-1}}, we obtain \textbf{Theorem \ref{mainth-2}}.

\textbf{Acknowledgment}
 Wu is supported by the Project of Stable Support for Youth Team in Basic Research Field, CAS(YSBR-001) and  the Fundamental Research Funds for the Central Universities.


\begin{thebibliography}{10}
\bibitem{Ca}E.Calabi,
\newblock ``Extremal K\"{a}hler metrics" in Seminar on
Differential Geometry.
\newblock {\em Ann. of Math. Stud. 102, Princeton Univ.
Press, Princeton},  259-290 (1982)

\bibitem{CL} W.H.Chen and  H.Z.Li,
\newblock Bonnet Surfaces and Isothermic Surfaces.
\newblock {\em Results. Math.}, 31, 40-52 (1997)

\bibitem{Ch0}X.X.Chen,
\newblock  Weak limits of Riemannian metrics in surfaces with integral curvature bound.
\newblock {\em  Calc. Var.}, 6, 189-226 (1998)

\bibitem{Ch1}X.X.Chen,
\newblock Extremal Hermitian metrics on Riemann surfaces.
\newblock {\em  Calc. Var. Partial Differential Equations 8} , no. 3, 191-232( 1999)

\bibitem{Ch2}X.X.Chen,
\newblock  Obstruction to the Existence of Metric whose
Curvature has Umbilical Hessian in a K-Surface.
\newblock {\em Comm. Anal. Geom.}, 8, no. 2, 267-299 (2000)

\bibitem{CCW}Q.Chen, X.X.Chen and Y.Y.Wu,
\newblock  The Structure of HCMU Metric in a K-Surface.
\newblock {\em Int. Math. Res. Not.}, 2005, no. 16, 941-958 (2005)

\bibitem{CW1}Q.Chen and Y.Y.Wu,
\newblock  Existences and Explicit Constructions of HCMU metrics on $S^2$ and $T^2$.
\newblock {\em  Pac. J. Math.}, 240, no. 2, 267-288 (2009)


\bibitem{CW2}Q.Chen and Y.Y. Wu,
\newblock  Character 1-form and the existence of an HCMU metric.
\newblock {\em Math. Ann.}, 351, no. 2, 327-345 (2011)

\bibitem{Xb}Q.Chen, Y.Y.Wu and B.Xu,
\newblock  On One-dimensional and singular calabi's extremal metrics whose gauss curvatures have nonzero umblical Hessians.
\newblock {\em Isr. J. Math.}, 208,  385-412 (2015)

\bibitem{LZ}C.S.Lin and X.H.Zhu,
\newblock Explicit construction of extremal Hermitian
metric with finite conical singularities on $S^2$.
\newblock {\em Comm. Anal. Geom.}, 10, no. 1, 177-216 (2002)

\bibitem{PW}C.K.Peng and Y.Y Wu,
 \newblock  A one-dimensional singular non-CSC extremal K$\ddot{a}$hler metric can be isometrically imbedded into $\mathbb{R}^{3}$ as a Weingarten surface.
 \newblock {\em Results. Math.}, 75, 133 (2020)


\bibitem{T}M.Troyanov,
 \newblock Prescribing curvature on compact surface with conical singularities.
\newblock {\em  Tran. Am. Math. Soc.}, 324(2),793-821 (1991)

\bibitem{WeiWu18}Z.Q.Wei and Y.Y.Wu,
\newblock One existence theorem for Non-CSC extremal K$\ddot{a}$hler metrics with singularities on $S^{2}$.
 \newblock {\em TJM} 22(1),55-62(2018)

 \bibitem{WeiWu22-1}Z.Q.Wei and Y.Y.Wu,
 \newblock Local isometric imbedding  of a compact Riemann surface with a singular  non-CSC extremal K$\ddot{a}$hler metric into  3-dimension space forms.
 \newblock {\em J. Geom. Anal.}, 32, 27 (2022)

\bibitem{WeiWu22-3}Z.Q.Wei and Y.Y.Wu,
\newblock  HCMU surfaces and Weingarten surfaces.
\newblock {\em J. Geom. Anal.}, 32, 199 (2022)

\bibitem{WeiWu22-4}Z.Q.Wei and Y.Y.Wu,
\newblock On isometric minimal immersion of a singular  non-CSC extremal K$\ddot{a}$hler metric into  3-dimension space forms.
 	arXiv:2112.15526

\bibitem{WZ}G.F.Wang and X.H.Zhu,
\newblock Extremal Hermitian metrics on Riemann surfaces with singularities.
\newblock {\em Duke Math. J.}, 104, 181-210 (2000)

\bibitem{MR} M.Dajczer and R. Tojeiro,
\newblock  Submanifold Theory: Beyond an Introduction.
 Springer, (2019)




\end{thebibliography}
\end{document}